\documentclass[12pt,leqno]{article}
\usepackage{amsmath,amssymb,amscd}
\numberwithin{equation}{section}

\newcommand{\Proof}{{\it Proof:\ }\ }

\renewcommand{\d}{\ensuremath{\mathrm{d}}}
\newcommand{\rond}{\circ}

\newtheorem{proposition}{Proposition}[section]
\newtheorem{lemma}{Lemma}[section]
\newtheorem{corollary}{Corollary}[section]
\newtheorem{definition}{Definition}[section]
\newtheorem{remark}{Remark}[section]
\newcommand{\ra}{\rightarrow}
\newcommand{\Hom}{\mathrm{Hom}}

\newcommand{\Id}{\mathrm{Id}}

\newcommand{\Tr}{\mathrm{Tr}}
\newcommand{\Spec}{\mathrm{Spec}}

\newcommand{\gr}{\mathrm{gr}}

\newcommand{\Z}{{\mathbb Z}}
\newcommand{\cqfd}{\ensuremath{\blacksquare}}

\newcommand{\A}{{\mathbb A}}
\newcommand{\X}{{\mathcal{X}}}

\renewcommand{\O}{\mathcal{O}}
\newcommand{\M}{\mathcal{M}}
\newcommand{\E}{\mathcal{E}}

\newcommand{\isom}{\stackrel{\sim}{\ra}}

\newcommand{\epi}{\ensuremath{\twoheadrightarrow}}

\mathcode`A="7041 \mathcode`B="7042 \mathcode`C="7043
\mathcode`D="7044 \mathcode`E="7045 \mathcode`F="7046
\mathcode`G="7047 \mathcode`H="7048 \mathcode`I="7049
\mathcode`J="704A \mathcode`K="704B \mathcode`L="704C
\mathcode`M="704D \mathcode`N="704E \mathcode`O="704F
\mathcode`P="7050 \mathcode`Q="7051 \mathcode`R="7052
\mathcode`S="7053 \mathcode`T="7054 \mathcode`U="7055
\mathcode`V="7056 \mathcode`W="7057 \mathcode`X="7058
\mathcode`Y="7059 \mathcode`Z="705A

\newcommand{\PP}{\mathcal{P}}

\newcommand{\Co}{\mathcal{C}}
\newcommand{\V}{\mathcal{V}}
\newcommand{\W}{\mathcal{W}}
\newcommand{\No}{\mathcal{N}ilp}
\newcommand{\dne}{{\mathcal{E}}nd}
\newcommand{\GG}{\mathbb{G}}
\newcommand{\Char}{\mathrm{Char}}

\begin{document}

\bibliographystyle{plain}
\title{\bf On the Hitchin morphism in positive characteristic}

\author{\textsc{Yves Laszlo} \and \textsc{Christian Pauly}}

\date{\today}

\maketitle

\begin{abstract}
Let $X$ be a smooth projective curve over a field of characteristic
$p>0$. We show that the Hitchin morphism, which associates to
a Higgs bundle its characteristic polynomial, has a non-trivial
deformation over the affine line. This deformation is
constructed by considering the moduli stack of $t$-connections
on vector bundles on $X$ and an analogue of the $p$-curvature,
and by observing that the associated characteristic polynomial 
is, in a suitable sense, a $p^{th}$-power.
\end{abstract}

\section{Introduction} Let $X$ be a smooth projective curve over an
algebraically closed field $k$ and let $\omega_X$ be its
canonical line bundle. The Hitchin morphism associates
to a rank $r$ vector bundle $E$ of degree zero and a Higgs
field $\phi: E \ra E \otimes \omega_X$ its characteristic polynomial,
denoted by $H(E,\phi)$, which lies in the affine space
$W = \oplus_{i=1}^r H^0(X, \omega_X^i)$. Thus one gets
a morphism
$$ H: \mathcal{H}iggs(r,X) \longrightarrow W$$
from the moduli stack of Higgs bundles to $W$, which becomes
universally closed, when restricted to the substack of semi-stable Higgs 
bundles \cite{Nit} \cite{Fal93}. Moreover, if
$k = \mathbb{C}$, it is shown \cite{hit87} that $H$
is an algebraically completely integrable system.

\bigskip
\noindent
In this note we show that the Hitchin morphism $H$ has a 
non-trivial deformation over the affine line if the 
characteristic of $k$ is $p >0$. More precisely, we
consider the moduli stack $\Co(r,X)$ of $t$-connections
$\nabla_t$ on rank $r$ vector bundles $E$ over $X$ with 
$t \in k$. A $t$-connection $\nabla_t$ can be thought of
as an ``interpolating'' object between a Higgs field
($t=0$) and a connection ($t=1$). Now one associates to
$\nabla_t$ a suitable analogue of the $p$-curvature
and it turns out (Proposition \ref{propdiv}) that its
characteristic polynomial is a $p^{\mathrm{th}}$-power.
This fact entails the existence of a morphism over $\A^1$
$$ H: \Co(r,X) \longrightarrow W \times \A^1$$
which restricts ($t=0$) to the Hitchin morphism for Higgs
bundles. Finally we prove that the restriction of $H$ to the
substack of semi-stable nilpotent $t$-connections is
universally closed. This provides a non-trivial  
deformation of the semi-stable locus of the global
nilpotent cone. This result can be considered as an
analogue of Simpson's result which says that the
moduli space of representations of the fundamental group has the
homotopy type of the global nilpotent cone.

\bigskip
\noindent
We thank G. Laumon for his interest and especially J. B. Bost  for
his help with the proof of Proposition \ref{propdiv}.

\section{Notations} We will denote by  $S$
a $k$-scheme with $k$ a field of characteristic $p>0$ and by $\X$
(resp. $X$) a smooth projective connected $S$-curve (resp.
$k$-curve).

\subsection{} \label{sect21} 
We denote by $f_S:  S\ra S$ the absolute Frobenius (which is
topologically the identity and the $p^{\mathrm{th}}$ power on
functions). We denote by $S^{(p)}$ the inverse image of $S$ by the
Frobenius $f_k$ of $k$. There exists a unique commutative diagram
\[\begin{array}{ccccccccc}
   S&\stackrel{F_S}{\ra}&S^{(p)}&\stackrel{\pi_S}{\ra}&S \\
   &\searrow&\downarrow&&\downarrow\\
   &&\Spec(k)&\stackrel{f_k}{\ra}&\Spec(k)
  \end{array}\]
  such that $\pi_S\rond F_S=f_S$. The $k$-morphism $F_S$ is called the
relative Frobenius.
If $S$ is defined by equations $f_j=\sum_I a_{I,j}x^I$ in an
affine space $\A^N$ with coordinates $(x_i)_{1\leq i \leq N}$ 
and $a_{I,j} \in k$, then $S^{(p)}$ is the subscheme of $\A^N$
defined by the equations $f_j^{[p]}=\sum_I
a_{I,j}^px^I$. In this case, $F_S: S\ra S^{(p)}$ is given by
$x_i \mapsto x_i^p$.

\subsection{} \label{sect22}
In the relative case $\X \ra S$, we denote by $\X^{(p)}$
the pullback of $\X$ by $f_S$. Assuming $S = \Spec(R)$, with 
$R$ a $k$-algebra, and $\X$ given by equations $f_j=\sum_I
a_{I,j}x^I$ in an affine space $\A^N_R$ with coordinates
$(x_i)_{1 \leq i \leq N}$ and $a_{I,j} \in R$, then 
$\X^{(p)}$ is the subscheme of $\A^N_R$ defined by the
equations $f_j^{[p]}=\sum_I a_{I,j}^px^I$. The $S$-morphism
$F_{\X/S}: \X \ra \X^{(p)}$ is given by $x_i \mapsto x_i^p$.

\subsection{}
 In the case $\X = X \times S$, the commutative diagram
\[
  \begin{array}{ccc}
    X^{(p)}\times S & \stackrel{(\pi_X,f_S)}{\ra}& X \times S \\
   \downarrow&&\downarrow\\
   S&\stackrel{f_S}{\ra}&S
  \end{array}
 \]
 defines a morphism $X^{(p)}\times S\ra (X\times S)^{(p)} = \X^{(p)}$, 
 which is an isomorphism thanks to the local description given in 
 (\ref{sect21}) and (\ref{sect22}). Therefore, we
 will identify the $S$-schemes $\X^{(p)}$ and $X^{(p)}\times S$ and 
 accordingly the relative Frobenius morphisms $F_{\X/S}$ and
 $F_X\times\Id_S$.

\section{The stack of bundles with connection}\label{stack:co}

Let $\X$ be a smooth projective curve over a basis $S$ of
characteristic $p$.

\subsection{} We denote by $\PP^1_{\X/S}$ the sheaf of principal parts of
degree $\leq 1$, namely the structural sheaf of the second order
infinitesimal neighborhood $\Delta^{(2)}$ of the diagonal in
$\X\times_S\X$ (\cite{eGAIV-4} section 16). The left
$\O_\X$-module structure is defined by the first projection and the
right one by the second projection. Let $1$ be the global section
of $\PP^1_{\X/S}$ defined by the constant function with value $1$
on $\Delta^{(2)}$. The restriction to the diagonal gives the exact
sequence
$$0\ra\omega_{\X/S}\ra\PP^1_{\X/S}\stackrel{\pi}{\ra}\O_{\X}\ra 0$$
of left $\O_{\X}$-modules. We tensorize this exact sequence with a 
vector bundle $\E$ on $\X$.
\begin{equation}\label{PE}
 0\ra\omega_{\X/S}\otimes\E\ra\PP^1_{\X/S}\otimes\E\ra\E\ra 0
\end{equation}
Note that $\PP^1_{\X/S}\otimes\E$ is the left-module tensor
product of the bimodule $\PP^1_{\X/S}$ by the left-module $\E$.
The following observation is classical (and tautological).

\begin{lemma}

(i) Let $\sigma\in\Hom(\E,\PP^1_{\X/S}\otimes\E)$ be a splitting of
(\ref{PE}). Then the morphism of sheaves $e\longmapsto
1.e-\sigma(e)$ takes its values in $\omega_{\X/S}\otimes\E$ and is
a connection $\nabla_\sigma$.

(ii) The map $\sigma\longmapsto\nabla_\sigma$ from splittings of
(\ref{PE}) to connections on $\E$ is bijective.
\end{lemma}\label{tauto}

 The functoriality of the sheaf of principal parts allows us to
 define the pull-back of a connection by a base change $S'\ra S$
 (which is a connection on the $S'$-curve $\X'=\X\times_S S'$).

\subsection{} Let $t \in H^0(S,\O)$ be a function on $S$.
We also denote by $t$ its pull-back to $\X$. 
We define the twisted sheaf of principal part $\PP^t_{\X/S}$
as the kernel of
\[
 \left\{ \begin{array}{ccc}
   \PP^1_{\X/S}\oplus \O_\X & \ra&\O_\X \\
    (p,e) & \longmapsto&\pi(p)-te\\
  \end{array}\right.
 \]
The sheaf $\PP^t_{\X/S}$ is an $\O_\X$-bimodule. The canonical
exact sequence
\begin{equation}\label{PtE}
 0\ra\omega_{\X/S}\otimes\E\ra\PP^t_{\X/S}\otimes\E\ra\E\ra 0
\end{equation}
is the pullback of (\ref{PE}) by the scalar endomorphism of $\E$
defined by $t$.

Recall that a $t$-connection on the vector bundle $\E$ on $\X$ 
is a morphism of $\O_S$-modules
$$\nabla_t:\ \E\ra \omega_{\X/S}\otimes \E$$
satisfying the (twisted) Leibniz rule
$$\nabla_t(fe)=t\d f\otimes e+f\nabla_t(e)$$
where $f$ and $e$ are local sections of $\O_\X$ and $\E$,
respectively. We say that the pair $(\E,\nabla_t)$ is a $t$-bundle. For
instance a $0$-connection is simply a Higgs field and, if $t$ is
invertible, $\nabla_t/t$ is a connection. Conversely, given a 
connection $\nabla$ on $\E$, then $\nabla_t := t \nabla$ is
a $t$-connection.

The twisted version of Lemma \ref{tauto} is

\begin{lemma}

(i) Let $\sigma\in\Hom(\E,\PP^t_{\X/S}\otimes\E)$ be a splitting of
(\ref{PtE}). Then the morphism of sheaves $e\longmapsto
1.e-\sigma(e)$ takes its values in $\omega_{\X/S}\otimes\E$ and is
a $t$-connection $\nabla_\sigma$.

(ii) The map $\sigma\longmapsto\nabla_\sigma$ from splittings of
(\ref{PtE}) to $t$-connections on $\E$ is bijective.
\end{lemma}\label{tauto-t}

\subsection{} Recall that $X$ is a smooth projective curve over $k$.
Let $S$ be a $t$-scheme,i.e., endowed with a global function $t$.
Let $\Co^*(S)$ be the category whose objects are pairs
$(\E,\nabla_t)$ where $\E$ is a degree $zero$ vector bundle on
$X_S=X\times S$ and $\nabla_t$ is a $t$-connection on $\E$ and
whose morphisms are isomorphisms commuting with the $t$-connections
(we still say flat morphisms). These categories with the obvious
inverse image functor define a fibred category over $Aff/\A^1$
which is obviously a stack, denoted by $\Co(r,X)$. Let $\M(r,X)$
denote the stack of rank $r$ vector bundles of degree zero over
$X$.

\begin{proposition} \label{prop31}
The forgetful morphism $\Co(r,X) \ra \M(r,X)\times\A^1$ is
representable.
\end{proposition}
\Proof Let $(\E,t)$ be an object of $\Hom_{\A^1}(S,\M(r,X)\times\A^1)$. Let
$\omega$ be the dualizing sheaf of $X_S\ra S$ (the pull-back of the
canonical bundle of $X$) and $\V$ be the dual of
$\mathcal{H}om(\E,\PP^t_{X_S/S}\otimes\E)\otimes\omega^{-1}$.

Let us consider the commutative diagram with cartesian square
\[
  \begin{array}{cccccc}
    &X_{S'}&\stackrel{g}{\ra}&X_S \\
    p'&\downarrow& &\downarrow&p\\
    &{S'}&\stackrel{f}{\ra}&S \\
    t'&\searrow& &\swarrow&t\\
    &&\A^1\\
  \end{array}
 \]

The category of triples $(\alpha,\E',\nabla_{t'})$ where $\alpha'$
is an isomorphism $\E'\isom g^*\E$ and $\nabla'$ a $t'$-connection
is equivalent with the subset (thought of as a discrete category)
of $$H(S')=\Hom(g^*\E,\PP^{t'}_{X_{S'}/S'}\otimes g^*\E)$$ mapping
to the identity of $ g^*\E$. The pull-back $\omega'$ of $\omega$ is
the dualizing sheaf of $X_{S'}$. Therefore
$$H(S')=H^0(p'_*\mathcal{H}om(g^*\V,\omega'))=\Hom(R^1p'_*g^*\V,\O_{S'})$$
by duality. By base change theory, $$R^1p'_*g^*\V=f^*R^1p_*\V$$
showing the canonical identification
$$H(S')=\Hom_S(S',V(R^1p_*\V)).$$ In the same way, $S'\ra \Hom(g^*\E,g^*\E)$
is representable by $V(R^1p_*\W)$ where $\W$ is the dual of
$\mathcal{H}om(\E,\E)\otimes\omega^{-1}$. The functorial morphism
induces an $S$-morphism
$$V(R^1p_*\V)\ra V(R^1p_*\W)$$
The identity of $\E$ defines a section of $V(R^1p_*\W)\ra S$ and
our $S$-category is represented by the fibred product
$$S\times_{V(R^1p_*\W)}V(R^1p_*\V).$$
\cqfd

Therefore, by general theory \cite{Lau-Mor99}, one gets
\begin{corollary}
The stack $\Co(r,X)$ is algebraic, locally of finite type over
$\A^1$.
\end{corollary}
\subsection{}\label{rappels:Psi(nabla)}

Let $\nabla$ be a connection on the rank $r$ vector bundle $\E$ over
$\X$. Recall that the $p$-curvature $\Psi(\nabla)$ of $\nabla$ is
the mapping (\cite{Kat70} section 5)
$$
  \Psi(\nabla):\left\{
  \begin{array}{ccc}
    T_{\X/S}&\ra&\dne_S(\E) \\
    D&\longmapsto&(\nabla(D))^p-\nabla(D^p)
  \end{array}\right.
$$
Recall that the relative tangent sheaf $T_{\X/S}$ is a sheaf of 
restricted $p$-Lie algebras.

\begin{remark}\rm
Notice that the curvature is the obstruction for the connection to
define a morphism of Lie algebra and that the $p$-curvature is the
obstruction for an integrable connection to define a morphism of
restricted $p$-Lie algebras.
\end{remark}

The additive morphism $\Psi(\nabla)$ is $\O_\X$-semi-linear,i.e.,
additive and $\Psi(\nabla)(gD) = g^p \Psi(\nabla)(D)$ for
$g$ and $D$ local sections of $\O_{\X}$ and $T_{\X/S}$ respectively
(\cite{Kat70} Proposition 5.2). Therefore it defines an
$\O_\X$-linear morphism
$$ \Psi(\nabla): T_{\X/S}\ra f_{\X*}\dne_S(\E)$$
where $f_\X$ is the absolute Frobenius of $\X$. We still denote by
$\Psi(\nabla)$ the corresponding section in 
$\Hom(\E,\E\otimes f_\X^*\omega_{\X/S})$ obtained by
adjunction.

\subsection{}\label{def:H} From now on we are interested 
in the case $\X=X\times S$. Let $q$ be the first projection. 
Since $f_{\X}^* \omega_{\X/S} = q^* \omega_X^p$, $\Psi(\nabla)$
can be thought of as a point of $\Hom(\E,\E\otimes q^*\omega_X^p)$. 
Let $\nabla_t$ be a $t$-connection on $\E$ with $t \in
H^0(S,\O)$. If $t\in \GG_m(S)$, then $\nabla_t/t$ is a 
connection on $\E$ and we have, for any local section
$D$ of $T_X$
$$t^p\Psi(\nabla_t/t)(D)= (\nabla_t(D))^p-t^{p-1}\nabla_t(D^p)
\in \dne_S(\E).$$

\begin{definition}
The $p$-curvature $\Psi_t(\nabla_t)$ of $\nabla_t$ is the additive
morphism
\[
  \Psi_t(\nabla_t):\
  \left\{
  \begin{array}{ccc}
    T_{\X/S} &\ra&\dne_S(\E)\\
    D &\longmapsto&(\nabla_t(D))^p-t^{p-1}\nabla_t(D^p)\\
  \end{array}
  \right.
 \]
\end{definition}

\begin{lemma}
The $p$-curvature of $\nabla_t$ is $\O_{\X}$-semi-linear.
\end{lemma}

\Proof We can adapt the proof of Proposition 5.2 \cite{Kat70} to 
the $t$-connection $\nabla_t$ taking into account that the
commutator of $\nabla_t(D)$ and $g$ in $\dne_S(\E)$
is $tD(g)$ (see formula (5.4.2) of \cite{Kat70}).
\cqfd

As in (\ref{rappels:Psi(nabla)}) we still denote by 
$\Psi_t(\nabla_t)$ the corresponding element in 
$\Hom(\E,\E\otimes q^*\omega_X^p)$.

\subsection{} Let us denote by $V$ the affine variety
$V=\oplus_{i=1}^rH^0(X,\omega_X^{pi})$ and by $V(S) = V \times S$.
\begin{definition}
Let $\nabla_t$ be a $t$-connection on $\E$. We denote by
$\Char(\nabla_t)$ the point of $V(S)$ defined by the coefficients
of the characteristic polynomial of the morphism $\Psi_t(\nabla_t):
\E \ra \E \otimes q^*\omega^p_X$.
\end{definition}

\begin{remark}\rm \label{rem32}
The functor $\nabla_t \longmapsto (\Char(\nabla_t),t)$ defines a
morphism of stacks
$$\underline{\Char}:\
\Co(r,X)\ra V\times\A^1$$

Given a Higgs bundle $(E,\nabla_0)$ over $X$, we observe that 
$$\Char(\nabla_0) = (H_0(\nabla_0))^p, $$
where $H_0$ is the Hitchin morphism as defined in \cite{hit87}.
As will be shown in the next section, $\Char(\nabla_t) \in V$
remains a $p^{\mathrm{th}}$ power when considering, more generally,
$t$-connections $\nabla_t$.
\end{remark}

\subsection{} Let us denote by $W$ the affine variety
$W = \oplus_{i=1}^rH^0(X,\omega_{X}^{i})$. The absolute
Frobenius morphism $f_X: X \ra X$
induces an injective $p$-linear map $f_X^*: W \hookrightarrow V$.

\begin{proposition}\label{propdiv}
There exists a unique morphism $H$ over $\A^1$
$$H:\ \Co(r,X)\ra
W \times\A^1$$ making the diagram
\[
  \begin{array}{ccccc}
    &\Co(r,X) \\
 & & & \\
    &\downarrow^H&\searrow^{\underline{\Char}}&\\
 & & & \\
&W\times\A^1&\hookrightarrow&V\times\A^1\\
  \end{array}
\]
commutative and such that the restriction $H_0$ of $H$ to the
fibre over $0 \in \A^1$ coincides with the Hitchin morphism.
\end{proposition}
\Proof
Let $(\E,\nabla_t)$ be a $t$-bundle over $X \times S$.
We can assume $S = \Spec(R)$, with $R$ integral since the 
stack $\Co(r,X)$ is smooth. 
To show the proposition, it will be enough to show that
$\Char(\nabla_t) \in W(S) \hookrightarrow V(S)$ or, equivalently,
that the global sections $\Char(\nabla_t) = (s_i) \in \oplus_i H^0(X,
\omega^{ip}_X) \otimes R$ descend to $X^{(p)} \times S$ 
by the relative Frobenius $F_X \times \Id: X \times S \ra X^{(p)}
\times S$. By Cartier's theorem (\cite{Kat70} Theorem 5.1) it 
suffices to check that for all $i$, $\nabla^{can} s_i = 0$, 
where $\nabla^{can}$ is the canonical connection on $\omega_X^{ip} =
F^*_X \omega^i_{X^{(p)}}$. The question is local on
$X$. Let $z$ be a local coordinate near a geometric point $x$ on
$X$. We choose trivializations of $\omega_X$ and $\E$. 
Then $\nabla^{can}$ is equal to the derivation $\d/\d z$, which we denote by
$\partial$, and the operator $\nabla_t(\partial) \in \dne(\E)$
can be written as 
\begin{equation} \label{exprcurv}
\nabla_t(\partial) = t \partial + A
\end{equation}
with $t \in R$ and $A$ an $r \times r$ matrix with entries in the
ring $R[[z]]$ (we work in a completion of the local ring $\O_{X,x}$).
We consider the $p$-curvature (note that $\partial^p = \partial$)
\begin{equation} \label{exprpcurv}
\Psi_t(\partial) = (t\partial + A)^p
\end{equation}
as an $r \times r$ matrix, which we simply denote by $\Psi_t$.
By $p$-linearity of the $p$-curvature $\Psi_t$ and by 
Cartier's theorem, one just has to show that the derivative
$$\partial.\det(\Id- T \Psi_t)= 0.$$

The following elegant argument is due to J.B. Bost. From \eqref{exprcurv}
and \eqref{exprpcurv} it is obvious that the commutator of
$\nabla_t(\partial)$ and $\Psi_t$ in $\dne(\E)$ is zero, hence
$$t [\partial,\Psi_t]=[\Psi_t,A]$$
Moreover one has $\partial. \Psi_t = [\partial,\Psi_t]$ where 
$\partial.\Psi_t$ is the matrix of derivatives of the
entries of $\Psi_t$. In particular, one gets for all $n \geq 0$
$$t \Tr(\Psi_t^n\partial.\Psi_t)= \Tr(\Psi^n_t [\Psi_t,A]) =
\Tr (\Psi_t^{n+1}A) - \Tr(\Psi_t^n A \Psi_t) = 0.$$
If $t=0$, then $\Psi_0 = A^p$ and we immediately see that
$\Tr(\Psi_0^n \partial.\Psi_0) = 0$. If $t \not= 0$, since
$R$ is integral, we get $\Tr(\Psi_t^n \partial.\Psi_t) = 0$.

Now we apply the formula $\partial. \det M = \det M \Tr(M^{-1} \partial.M)$
to the matrix $M = \Id - T \Psi_t$ and we write $M^{-1} = \sum_{n \geq 0}
T^n \Psi_t^n$, 
$$\partial.\det(\Id -T\Psi_t) = -T \det(\Id -T \Psi_t)
\sum_{n \geq 0} T^n \Tr(\Psi_t^n \partial.\Psi_t).$$
Since, for all $t$ and $n$, the elements $\Tr(\Psi_t^n \partial.\Psi_t)$
are zero, we get the result. $\blacksquare$

\section{The stack of nilpotent connections}
We consider the embedding $\A^1 \hookrightarrow W \times \A^1$ by the
zero section and we denote by $\No(r,X)$ the corresponding fibre 
product $\Co(r,X)\times_{W \times\A^1}\A^1$,
$$\begin{CD}
\No(r,X) @>>> \Co(r,X) \\
@VVHV  @VVHV \\
\A^1 @>0>> W \times \A^1 
\end{CD}$$
The category $\No(r,X)(S)$ over a $t$-scheme
$S$ is the category of $t$-bundles $(\E,\nabla_t)$ over $X \times S$
such that
$(\Psi_t(\nabla_t))^r\in \Hom(\E,\E\otimes q^*\omega_X^{pr})$ is
zero. We say that $\nabla_t$ is nilpotent. The stack of
nilpotent connections is
far from being reduced if $r>1$. We define the exponent of
nilpotence of the $t$-connection $\nabla_t$ as the smallest
integer $l \geq 1$ such that  $(\Psi_t(\nabla_t))^l = 0$ and
we denote by $\No_l(r,X)$ the substack of $\No(r,X)$ corresponding
to nilpotent $t$-connections of exponent $\leq l$. Thus we
get a filtration
$$\No_1(r,X) \subset \cdots \subset \No_r(r,X) = \No(r,X)$$
and by Theorem 5.1 \cite{Kat70} we have an isomorphism induced by
the relative Frobenius $F_X$
\begin{equation}\label{pbrf}
\M(r,X^{(p)}) \isom \No_1(r,X), \qquad E \mapsto 
(F^*_X E, \nabla^{can}). 
\end{equation}

\subsection{} Following \cite{Sim94}, we say that $(\E,\nabla_t)$
over $X \times S$ (or simply $\nabla_t)$ is semi-stable 
if for all geometric points
$\bar s$ of $s \in S$, any proper $\nabla_{t,\bar s}$-invariant
subsheaf of $\E_s$ has negative degree (recall that $\E$ is of degree
$0$). We denote by $\Co^{ss}(r,X)$ (resp. $\No^{ss}_l(r,X))$ the substack
of $\Co(r,X)$ (resp. $\No_l(r,X)$) parameterizing semi-stable
$t$-connections (resp. semi-stable $t$-connections which are
nilpotent of exponent $\leq l$).

\begin{lemma}
The natural morphism $$\Co^{ss}(r,X)\ra\Co(r,X)$$ is an open
immersion.
\end{lemma}
\Proof Let $(\E,\nabla_t)$ be an $S$-point of $\Co(r,X)$. One has
to prove that the semi-stable sublocus in $S$ is open. Because
$\Co(r,X)$ is locally of finite type, one can assume that $S$ is of
finite type. By the standard argument for boundedness, 
there exists $n>>0$ such
that every subsheaf of $\E_s, s\in S(k)$, has slope $\leq n$. Let
$I$ be the finite set
$$I=\{(r,d)\in\{1,\ldots,r-1\}\times\Z\ \mathrm{such\ that}\ 
0<\frac{d}{r}\leq n\}.$$
Let $\mathcal{Q}^I$ be the relative Hilbert scheme parameterizing
subsheaves $F_s$of $\E_s$ with Hilbert polynomial in $I$. Because
$I$ is finite, $\mathcal{Q}^I$ is proper over $S$. But
$F_s$ is flat (i.e. $\nabla_{t,s}$-invariant) if and only if
the linear morphism $F_s\ra(\E_s/F_s)\otimes\omega$ is zero. As 
in the proof of Proposition \ref{prop31},
it follows that
the sublocus of $\mathcal{Q}^I$ parameterizing flat subsheaves is
closed. Its image in $S$ is closed and is exactly the 
locus where $(\E_s,\nabla_{t,s})$
is not semi-stable.\cqfd
\subsection{} Let $E$ be a vector bundle on $X$ with
Harder-Narasimhan filtration
$$0=E_0\subsetneqq E_1\subsetneqq\ldots\subsetneqq E_l=E.$$
We denote by $\mu_+(E)$ (resp. $\mu_-(E)$) the slope of the
semi-stable bundle $E_1/E_0$ (resp. $E_l/E_{l-1}$). One has by
construction
$$\mu_-(E)\leq\mu(E)\leq\mu_+(E)$$
with equality if and only if $l=1$ and $E$ semi-stable.
\begin{lemma}
The algebraic stack $\Co^{ss}(r,X)$ is of finite type.
\end{lemma}
\Proof By the boundedness of the family of semi-stable bundles of
given slope and rank, it is enough to find a bound for
$\mu_+(E)-\mu_-(E)$ where $(E,\nabla_t)$ is semi-stable on $X$.
Assume $l>1$ and let us consider the linear morphism
$$\bar\nabla_i:\ E_i\ra (E/E_i)\otimes\omega_X$$
induced by $\nabla_t$ for $1<i<l$. By construction,
$\mu(E_i)>\mu(E)=0$. If $\bar\nabla_i$ vanishes, $E_i$ is flat and
therefore of non positive slope, a contradiction. It follows that
$$\mu(E_i/E_{i-1})=\mu_-(E_i)\leq\mu_+((E/E_i)\otimes\omega_X)
=\mu(E_{i+1}/E_{i})+2g-2.$$
In particular, one has
$$0\leq\mu_+(E)-\mu_-(E)\leq(l-1)(2g-2) \leq (r-1)(2g-2).$$\cqfd

\section{} We give a properness property
of the deformation of the nilpotent cone of Higgs bundles induced by $H$. The
precise statement is
\begin{proposition}
For all $l \geq 1$, the morphism of stacks of finite type 
$$H:\ \No^{ss}_l(r,X) \ra \A^1$$ 
is universally closed.
\end{proposition}

\begin{remark} \rm
We recall that, if $l=1$, the fibres $H^{-1}(0)$ and $H^{-1}(1)$
are equal to the moduli stacks $\M^{ss}(r,X)$ and $\M^{ss}(r,X^{(p)})$
of semi-stable vector bundles, respectively (see \eqref{pbrf}).
\end{remark}

By \cite{Lau-Mor99}, it is enough to prove the following
semi-stable reduction theorem.

\begin{proposition}
Let $R$ be a complete discrete valuation ring with field of
fractions $K$ and $t\in R$. Let $\X$ be a smooth, projective,
connected curve over $R$ and $(\E_K,\nabla_K)$ a semi-stable
$t$-bundle on $\X_K$ with nilpotent $p$-curvature of exponent $l$. 
Then there exists a semi-stable $t$-bundle $(\E,\nabla)$ with nilpotent
$p$-curvature on $\X$ extending $(\E_K,\nabla_K)$ and of exponent $l$.
\end{proposition}\label{th:propre}

{\it Proof:} We adapt Faltings's proof, based on ideas of Langton,
of the properness of the moduli of Higgs bundles (\cite{Fal93} Theorem
I.3) to $t$-bundles. If $t = 0$, the result is a particular case
of Faltings's theorem. So we assume $t\not=0$. We denote by $F:\
\X\ra\X^{(p)}$ the relative Frobenius, by $F_K:\ \X_K\ra\X_K^{(p)}$
its restriction to the generic fiber, and by $\Psi_K$ the 
$p$-curvature of the connection $\nabla_K/t$ on $\E_K$. By assumption
$\Psi_K^l=0$.

\begin{lemma} \label{lemma51}
Suppose that $(\E_K,\nabla_K)$ is a (not necessarily semi-stable)
$t$-bundle of exponent of nilpotence $l$ on $\X_K$. Then it can be
extended to the whole $\X$ with exponent $l$.
\end{lemma}

{\it Proof of the lemma:} We proceed by induction on $l$.
If $l=1$, $\Psi_K = 0$ and, by Theorem 5.1 \cite{Kat70},
there exists a vector bundle $\E^{(p)}_K$ on $\X^{(p)}_K$ such that
$F_K^*\E^{(p)}_K\isom\E_K$ and $\nabla_K/t$ becomes the canonical
connection $\nabla^{can}$ on $F_K^*\E^{(p)}_K$. Now 
$\E^{(p)}_K$ extends to a vector bundle $\E^{(p)}$ on $\X^{(p)}$
and the $t$-bundle $(F^*\E^{(p)},t\nabla^{can})$ extends $(\E_K,\nabla_K)$,
since $t\in R$.

Now assume that $l>1$ and that the lemma holds for $t$-bundles
with exponent of nilpotence $<l$. Let $F^n_K$ be the
kernel of $\Psi_K^n$ for $n=1,\ldots,n$. Thus we have a filtration
$$ 0 = F^0_K \subset F^1_K \subset \ldots \subset F^{l-1}_K
\subset F_K^l = \E_K.$$
Since $\Psi_K$ and $\nabla_K/t$ commute (Proposition 5.2 \cite{Kat70}),
the connection
$\nabla_K/t$ maps $F^n_K$ to $F^n_K\otimes\omega_{\X_K/K}$.
By construction,
the induced connection $\nabla_K^F$ on $F_K^{l-1}$ is of 
exponent $l-1$. Moreover $\nabla_K/t$ induces
a connection $\nabla_K^{\gr}$ on the quotient $\gr_K := \E_K/F_K^{l-1}$ 
, which is of exponent $1$, since $\Psi_K$ maps $\E_K$ to
$F_K^{l-1}$.

By induction the $t$-bundles $(F^{l-1}_K,\nabla_K^F)$ and
$(\gr_K,\nabla_K^{\gr})$ have models on $\X$, which we denote by
$(F^{l-1},\nabla^{F})$ and $(\gr,\nabla^\gr)$, respectively.

Let $x$ be the generic point of the special fiber of $\X$. The
$R$-algebra $A=\O_{\X,x}$ is a discrete valuation ring with field
of fraction $L=k(\eta)$ where $\eta$ is the generic point of $\X$.
By Proposition 6 \cite{Lan75}, it is enough to find an $A$-lattice $E_x$ in the
$r$-dimensional $L$-vector space $\E_{K,\eta}$ such that
\begin{equation} \label{exisext}
\nabla_K/t(E_x)\subset E_x\otimes_A\omega_{x}
\end{equation}
where $\omega_{x}=\omega_{\X/R,x}$. 

Now, $\E_K$ is an extension of $\gr_K$ by $F^{l-1}_K$ whose
extension class $e_K$ lives in
$$H^1(\X_K,Hom(\gr_K,F^{l-1}_K))= H^1(\X,Hom(\gr,F^{l-1}))\otimes
K.$$
Let $\pi$ be a uniformizing parameter of $R$. If $m$ is large enough,
$\pi^me_K$ defines an extension
\begin{equation}\label{ext:grF}
  0\ra F^{l-1}\ra\E\ra\gr^l\ra 0
\end{equation}
on $\X$ together with an isomorphism
$$\E\otimes K\isom\E_K.$$
The localization of (\ref{ext:grF}) at $x$ is split ($A$ is
principal) and, with respect to this splittings, $\nabla_K/t$ is of
the form
\[
\nabla_K/t=\left(
\begin{array}{ccc}
\nabla^{F}&0\\
0&\nabla^\gr\\
\end{array}
\right)
+
\left(\begin{array}{ccc}
0&M\\
0&0\\
\end{array}
\right)
\]
where $M\in L\otimes\Hom(\gr_x,F^{l-1}_x\otimes\omega_x)$. Let $m$
be an integer such that $\pi^mM\in
\Hom(\gr_x,F^{l-1}_x\otimes\omega_x)$ and let
$\E_x^m\subset\E_{K,\eta}$ be the lattice defined by the pull-back
of
$$0\ra F^{l-1}_x\ra\E_x\ra\gr_x\ra 0$$ by
\[\left\{\begin{array}{ccc}
\gr_x&\ra&\gr_x\\
g&\longmapsto&\pi^mg\\
\end{array}
\right.
\]
Hence, by construction, the $A$-lattice $\E^m_x$ satisfies
\eqref{exisext}. Then the vector bundle $\E^m$ defined by $\E_K$ and 
$\E_x^m$, and the $t$-connection $\nabla_K=t\nabla_K/t$ define
a $t$-bundle extending $(\E_K,\nabla_K)$ of exponent $l$. 
\cqfd

\bigskip

The rest of the proof goes exactly as in \cite{Fal93} (or \cite{Lan75}). 
Let us for the
convenience of the reader recall the argument. Let $s$ be the
closed point of $\Spec(R)$ and $\bar s$ a geometric point over it.
As usual every $t_{\bar s}$-bundle over $\X_{\bar s}$ has a unique
flat (i.e., invariant) subbundle of maximal slope. From this
uniqueness, it follows easily (adapt the arguments of \cite{Ses82}
where the case $t_s=0$ is treated, for instance), this bundle is in
fact defined over $\X_s$. In particular, there is a notion of
Harder-Narasimhan filtration of $t_s$-bundles. The set of slopes of
the first term $F^\alpha(\E_s,\nabla_s)$ of the Harder-Narasimhan
filtration of the $t_s$-bundle $(\E_s,\nabla_s)$ where
$(\E,\nabla)$ runs over the models of $(\E_K,\nabla_K)$ with
$p$-curvature of exponent $l$,
whose existence is proved in Lemma \ref{lemma51}, has
certainly a smallest element $\alpha$. Among these models of
minimal slopes $\alpha$, pick-up one, $(\E^0,\nabla^0)$ say, such
that $F^\alpha=F^\alpha\E^0_s$ has minimal rank.
We will show that the $t$-bundle $(\E^0, \nabla^0)$ is
a semi-stable model. Suppose that the contrary holds,i.e.,
$\mathrm{rk}(F^\alpha) < r$.

The connection $\nabla^0$ induces a connection
$\nabla^1$ on the kernel $\E^1$ of the surjection
$\E^0\epi\E^0_s/F^\alpha$ and therefore we get a new model
$(\E^1,\nabla^1)$ with $\E^1$ locally free. Notice that 
$\nabla^1$ still has $p$-curvature of exponent $l$. Because
$${\mathcal{T}}or_1^{\O_\X}(\E^0_s/F^\alpha,\O_{\X_s})=
\E^0_s/F^\alpha,$$
the restriction of the exact sequence
$$0\ra\E^1\ra\E^0\ra\E^0_s/F^\alpha\ra 0$$
to $\X_s$ gives an exact sequence of flat bundles
$$0\ra\E^0_s/F^\alpha\ra\E^1_s\ra F^\alpha\ra 0.$$

Let $F^\beta_1$ be the first term of the Harder-Narashiman
filtration of $(\E^1_s,\nabla_{t,s})$. By minimality of $\alpha$,
one has $\alpha\leq\beta$. Observe that the Harder-Narashiman
filtration $\E^0_s/F^\alpha$ has semi-stable subquotients of slopes
$<\alpha\leq\beta$. Therefore, every morphism from a stable bundle
of slope $\beta$ to $\E^0_s/F^\alpha$ is zero. It follows
immediately that $\beta=\alpha$ and that $F^\beta_1$ injects to
$F^\alpha$. By minimality of the rank of $F^\alpha$, this injection
is an isomorphism. Therefore $F^\alpha$ is a subsheaf of $\E^1$. By
induction, one can construct the subsheaf $\E^{n+1}$ of $\E^0$ as
the kernel of $\E^n\ra\E^n_s/F^\alpha$ with a
$(t\mod\pi^{n+1})$-connection $\nabla^{n+1}$ induced by $\nabla^n$.
The point is that one has the exact sequence (of subsheaves of
$\E^0$)
$$0\ra \pi\E^n\ra\E^{n+1}\ra F^\alpha\ra 0.$$
In particular, $\pi^n\E^0$ is contained in $\E^n$ and
$\E^n/\pi^n\E^0$ injects in $\E^0/\pi^n\E^0$. Taking the direct
limit, we get a subsheaf of the completion of $\E^0$ along $\X_s$
with restricts to $F^\alpha$ on the special fiber. It comes with a
$t$-connection induced by $\nabla^0$, and both the sheaf (which is
a vector bundle) and the $t$-connection algebraizes by properness
of $\X\ra\Spec(R)$. Moreover, this subsheaf is certainly
$\nabla$-invariant. By flatness over $R$, the slope of its generic
fiber is $\alpha >0$, contradicting $(\E_K,\nabla_K)$
semi-stable.
\cqfd

\bigskip

\flushleft{Yves Laszlo \\
Universit\'e Paris-Sud \\
Math\'ematiques B\^atiment 425 \\
91405 Orsay Cedex France \\
e-mail: Yves.Laszlo@math.u-psud.fr}

\bigskip

\flushleft{Christian Pauly \\
Laboratoire J.-A. Dieudonn\'e \\
Universit\'e de Nice Sophia Antipolis \\
Parc Valrose \\
06108 Nice Cedex 02 France \\
e-mail: pauly@math.unice.fr }

\end{document}